\documentclass[letterpaper, 10 pt, conference]{ieeeconf}  

\IEEEoverridecommandlockouts                              
\usepackage{amsmath} 
\usepackage{amssymb}  

\usepackage{color}
\usepackage{soul}
\usepackage{graphicx}
\usepackage{caption}
\usepackage{subcaption}

\title{\LARGE \bf
Evaluation of Communication Issues in Primal-Dual-Based Distributed Energy Resource Management Systems (DERMS)}

\author{
Joshua Comden, Jing Wang, and Andrey Bernstein
\thanks{The authors are with the Power Systems Engineering Center at the National Renewable Energy Laboratory, Golden, CO, USA. Emails: \{Joshua.Comden, Jing.Wang, Andrey.Bernstein\}@nrel.gov}
}

\begin{document}

\maketitle
\thispagestyle{empty}
\pagestyle{empty}

\begin{abstract}
With the increasing adoption of distributed energy resources (DERs) in distribution networks, distributed energy resource management systems (DERMS) are becoming an attractive option to coordinate the control of DERs, especially primal-dual-based DERMS, which is a well-developed class.
To help reduce the uncertainty in commercializing DERMS, we evaluate and find the upper limits of severity for some common issues that can occur within the communication systems upon which a primal-dual-based DERMS relies; these issues include packet losses, link failures, and delays.
This information can help a DERMS operator decide on the specifications for the communication infrastructure it will invest in for the DERMS.
Our evaluation is based on numerical simulations of a real-world feeder with approximately 2,000 nodes and a DERMS that controls 163 photovoltaic generators and 140 energy storage batteries.
In general, we find that the DERMS becomes less resistant to communication issues as the information flows closer to the DERs implementing the controllable power injections.
\end{abstract}

\section{Introduction}
\label{sec:intro}
The increasing number of distributed energy resources (DERs) being connected to distribution networks is necessitating their coordinated control; otherwise, allowing DERs to remain uncoordinated can cause voltage or power line capacity violations within distribution networks.
A distributed energy resource management system (DERMS) is a general control framework that can be used by a utility or other third party to not only coordinate the control of DERs but also provide services to the grid, such as voltage support.

A well-developed class of DERMS, characterized by its hierarchical primal-dual control structure, uses a centralized coordinator to collect measurements associated with grid services from various locations on the distribution feeder and to send control signals to its distributed local controllers, each one determining the set points for the DER(s) it controls.
The DERMS has been designed to provide various combinations of grid services, including voltage support~\cite{dall2016optimal,wang2021performance,padullaparti2022evaluation}; voltage support and virtual power plants (VPPs)~\cite{dall2017optimal,gan2021cyber,wang2020performance1,wang2021voltage};
and voltage support, line current constraints, and VPPs~\cite{bernstein2019real}.
A primal-dual-based DERMS heavily relies on communication infrastructure to collect grid service measurements, send control signals, and implement DER set points.

When a system operator decides to invest in a DERMS, it needs to balance the cost of the communication infrastructure with the quality of the grid services that it intends to provide;
however, communication systems can experience a variety of issues---including packet losses, link failures, and delays---which can degrade the performance of the control systems built upon them, such as DERMS~\cite{gan2021cyber}.
Therefore, our goal is to find the severity limits of these communication issues with respect to the functionality of primal-dual-based DERMS so that an operator can choose the cost-appropriate communication infrastructure that respects these limits.
This reduces the investment uncertainty for an operator and thus promotes the commercialization of DERMS.

Although most previous evaluations of primal-dual-based DERMS assume that communication systems are completely reliable and have no significant delays, some more recent work has started to investigate the effects that more realistic communication system properties have on DERMS;
however, they are more focused on average cases rather than finding the upper limits on functionality.
For example, \cite{gan2021cyber} studies communication delays on the order of milliseconds instead of the order of seconds or tens of seconds, as we do in this paper.
Other works look at the opposite impact perspective to us by focusing on how different properties of DERMS affect the operational quality of the communication systems.
The study done by \cite{zhang2019analysis} investigates how DER penetration affects the size of the delays and packet loss rate.

To help a primal-dual-based DERMS operator make the appropriate communication infrastructure decisions with respect to communication issues---including the upper bounds on packet losses, link failures, and delays---this paper makes the following contributions:
\begin{enumerate}
    \item We describe in detail the communication architecture of primal-dual-based DERMS and categorize the communication channels (Section \ref{sec:communication}).
    \item We describe and model the communication issues so that they can be simulated in our numerical evaluation (Section \ref{sec:issues}).
    \item We describe the evaluation metrics used to determine when the DERMS is functional or not (Section \ref{sec:metrics}).
    \item We use numerical simulations of a real-world feeder and the system's communication channels to find the upper bounds on the properties of the communication issues (Section \ref{sec:perf_eval}).
\end{enumerate}

\section{Communication Architecture}
\label{sec:communication}
\begin{figure}
    \centering
    \includegraphics[width=0.98\columnwidth]{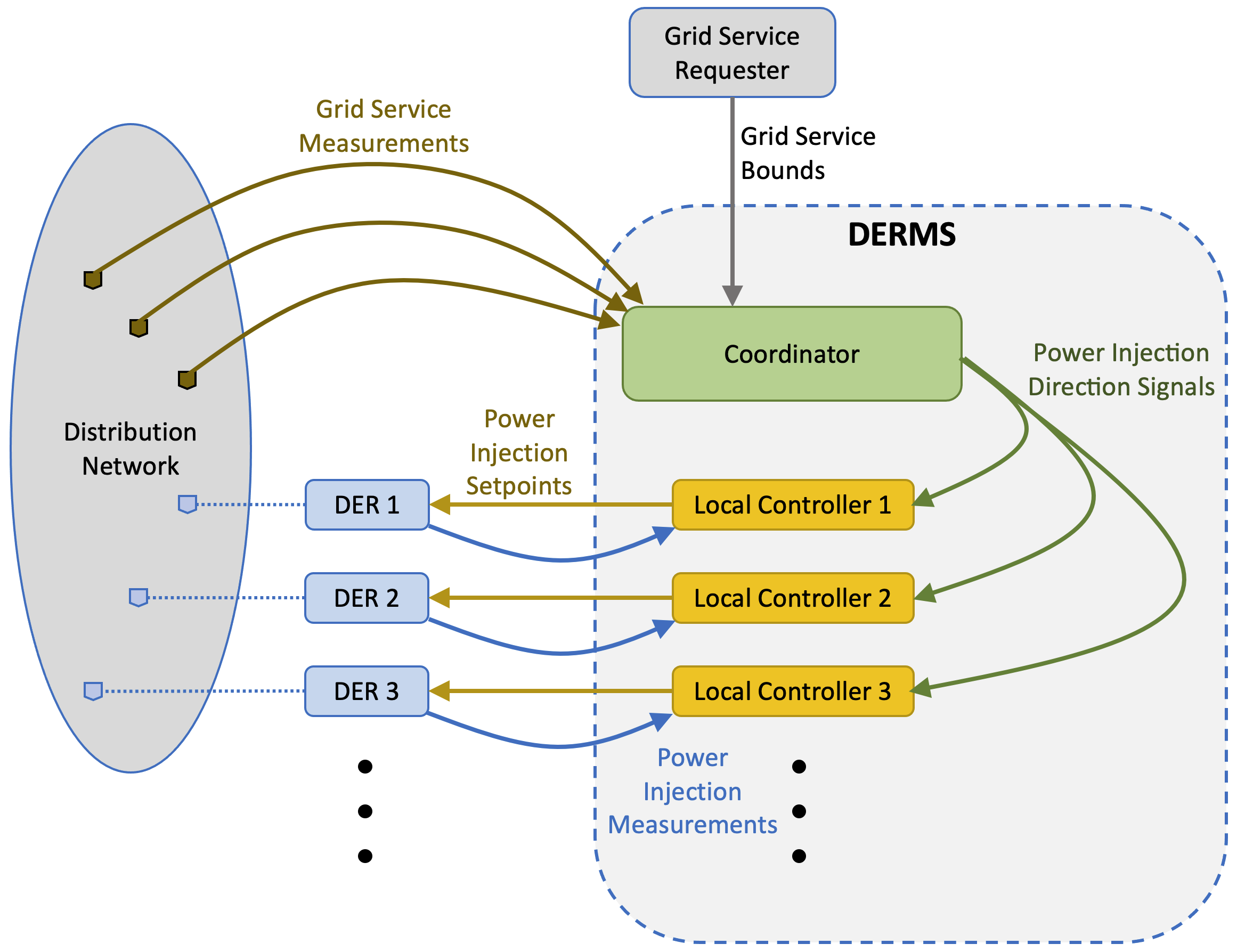}
    \caption{Information flow of primal-dual-based DERMS}
    \label{fig:DERMS_diagram}
\end{figure}

A primal-dual-based DERMS works as a hierarchical feedback controller~\cite{bernstein2019online} of a distribution network where the power injections supplied by the DERs are the control variables, and the grid service measurements (e.g., voltage measurements) are the feedback variables to be controlled.
The DERMS operator supplies grid services to an entity that sets the time-varying bounds for the grid service measurements.
Note that there are other primal-dual DERMS structures that are fully decentralized that do not require a coordinator (e.g., \cite{liu2018hybrid,patari2022distributed}), but they require heavier use of the communication infrastructure.

We classify the DERMS communication channels by the information being sent across them into the following categories, as shown in Figure \ref{fig:DERMS_diagram}:
\begin{enumerate}
    \item Grid service measurements
    \item Power injection direction signals
    \item Power injection measurements and set points.
\end{enumerate}
The first category uses a neighborhood area network (NAN) to send grid service measurements from across the distribution network to the DERMS coordinator, which converts the grid service measurements that are outside their associated bounds into control signals that are individualized power injection direction signals.
The second category uses a NAN to send the power injection direction signals from the coordinator to the distributed local controllers, which update the power injection set points of the DERs.
Finally, the third category uses a home area network (HAN) to send the power injection set points from the local controller to the DER(s) it controls and the power injection measurements from the DER(s) to its local controller(s).

\section{Communication Issue Models}
\label{sec:issues}
In this section, we describe the three communication issues that we study and how we model them in the numerical simulations.
For each communication issue, we evaluate each communication channel category separately by only applying the issue to that individual category and observing the impacts.

\subsection{Packet Losses}
When using communication networks to pass messages between different components of the DERMS, messages can occasionally drop and not make it to their destination.
These usually happen without warning and are one-off events;
thus, we model this phenomenon as a Bernoulli process, where every time a message is sent, there is a specific probability that the message will not be delivered to its destination~\cite{liu2018hybrid}.

For each category, we evaluate the upper bound on the probability of a packet loss that will keep the DERMS functional. 
We vary the number of channels in the category that are set to have a random pack loss to find the upper bound on the probability as a function of the number of channels.
All other channels are assumed to have no random packet loss.

\subsection{Link Failures}
In communication networks, links can randomly fail for several minutes at a time.
We simulate this phenomenon similarly to the random packet drop simulations.
If a communication link is functional, after it is used to send a message, a Bernoulli random variable is used to decide whether the link will become nonfunctional.
If it becomes nonfunctional, then the communication channel cannot be used to send messages for a set amount of time.

For each channel category, we evaluate the upper bound on the probability as a function of the duration of time that the link will be nonfunctional, assuming every link in the channel category can become nonfunctional.

\subsection{Delays}
When sending messages across a communication network, there is an interval of time from when the data in the message were created to when they are used in the next control component, called a delay or latency.  This can have adverse oscillatory effects on controlling systems due to the feedback from delayed actions.

Although other evaluations use a probabilistic model to study the effects of delays (e.g., \cite{gan2021cyber}), we use a constant delay size for all of the communication channels within the category we evaluate (e.g., \cite{magnusson2020distributed}); setting all channels within a category to be at the upper bound is the worst-case scenario.
For each channel category, we evaluate the upper bound on the delay that can keep the DERMS functional by applying the delay to all channels within the category.

\section{Communication Evaluation Metrics}
\label{sec:metrics}
As we increase the severity of the communication issues to find their limits for the different communication channel categories, we need to specify an evaluation metric that determines whether the DERMS is operating in a functional manner.

After a large disturbance, a functioning DERMS should be able to bring the grid service measurements back to within or at least near their bounds within a reasonable amount of time;
thus, we define a metric based on this condition.
The metric is specified by a time duration and a tolerance for each grid service measurement.
We consider a DERMS to be functional if within the time duration after a disturbance, it is able to get all of its grid service measurements to be within their bounds that are widened by the tolerance.
The specific values of the time duration and the grid service measurement tolerances could be determined by the grid service contract that the DERMS is operating under.

\section{Numerical Evaluation}
\label{sec:perf_eval}
In this section, we use the metrics described in Section \ref{sec:metrics} to find the severity limits of the communication issues described in Section \ref{sec:issues} for the communication channel categories described in Section \ref{sec:communication} with numerical simulations of a real-world feeder and a primal-dual-based DERMS.

\subsection{Setup}

\begin{figure}
    \centering
    \includegraphics[width=0.98\columnwidth]{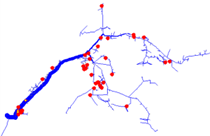}
    \caption{Feeder topology with the controllable DERs in red.}
    \label{fig:topology}
\end{figure}

The control system to be evaluated is a primal-dual-based DERMS that is contracted to provide voltage support and VPP grid services for a three-phase distribution feeder.
The voltage support lower and upper bounds are set at a constant 0.95 p.u. and 1.03 p.u., respectively.
The VPP is defined as $\pm$0.01 MW bounds from a time-varying set point for the feeder head powers.
The feeder contains approximately 2,0000 nodes, with an aggregate peak of 4.6 MW, and it is modeled on a utility system in Colorado.
Each node containing a DER(s) has a DERMS local controller, and each DER is equipped with a smart inverter to inject real and reactive power.  

The DERMS controls 163 curtailable photovoltaic (PV) solar generators (see Figure \ref{fig:topology}) with capacities that range from 0.04 kW to 34 kW, averaging 10 kW, and that double as the sensors sending the voltage magnitude measurements back to the coordinator; 140 of them have an individually controllable energy storage battery, with storage capacities that range from 13.5 kWh to 54 kWh, averaging 19 kWh.
The local DER costs are modeled by quadratic functions that represent the cost of curtailing away from its available power in the case of a PV generator and the cost of diverging the state of charge from its preferred value in the case of a battery.
The quadratic cost coefficients for each DER are scaled by its capacity.
More specifics on the cost models can be found in \cite{padullaparti2022evaluation,bernstein2019real}.

\begin{figure}
    \centering
    \includegraphics[width=0.98\columnwidth]{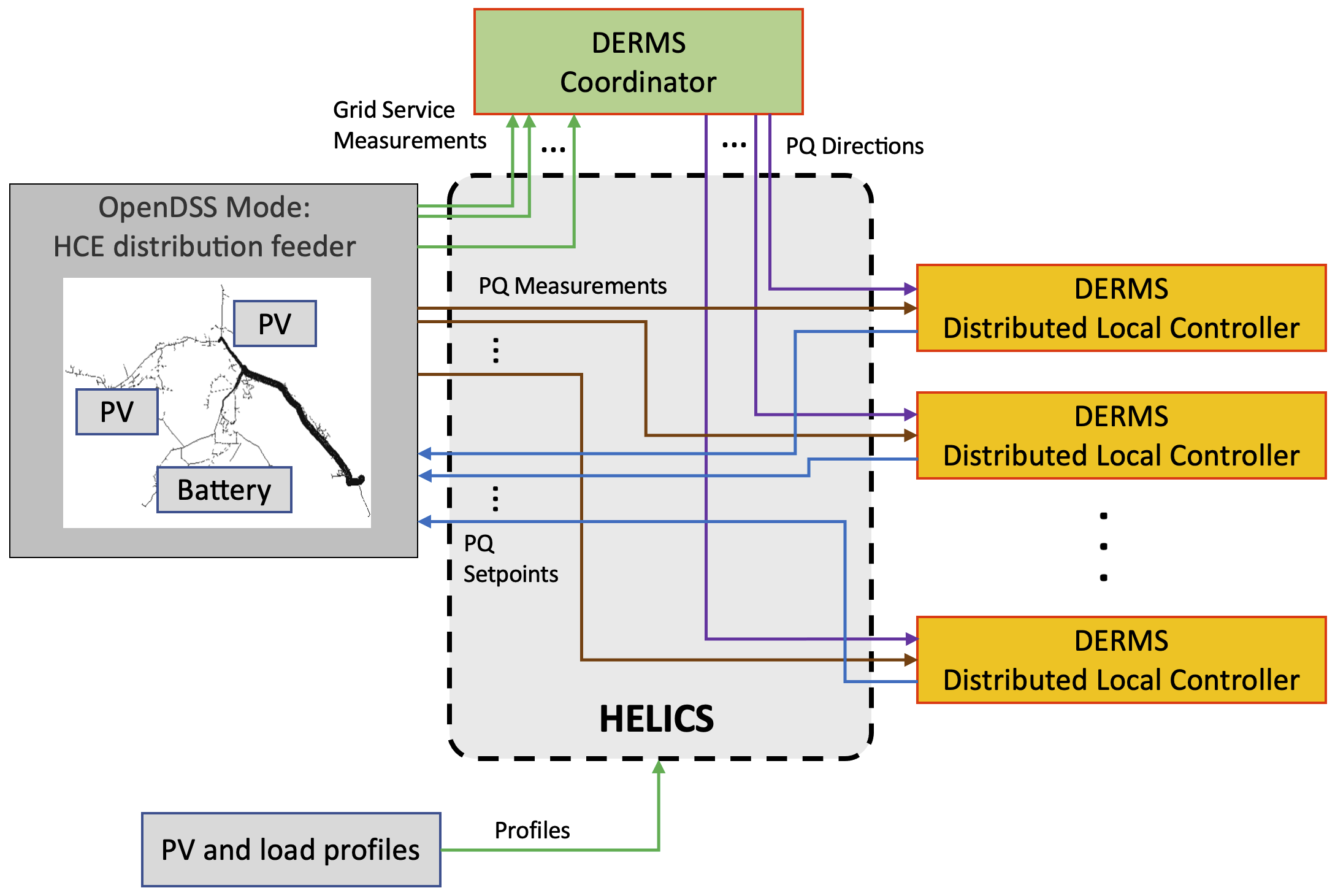}
    \caption{Co-simulation diagram with HELICS.}
    \label{fig:HELICS_diagram}
\end{figure}

The feeder, the coordinator, the distributed local controllers, and the DERs are co-simulated with the Hierarchical Engine for Large-scale Infrastructure Co-Simulation (HELICS)~\cite{palmintier2017design}, which coordinates the information sent between each module and the execution time.
Figure \ref{fig:HELICS_diagram} shows the flow of information between the modules with HELICS.
The feeder and DERs are simulated in OpenDSS every 2 seconds as a quasi-steady-state time series.
The utility provided the load data at a 15-minute time resolution and the PV generation data at a 1-minute time resolution from their advanced metering infrastructure.
The DERMS coordinator updates its internal (dual) variables every 1 minute, and each local controller updates its (primal) variables every time it receives a power injection measurement from the DER(s) it controls.
Unless otherwise noted, the grid service measurements (voltage magnitudes and feeder head powers) are sent to the coordinator every 1 minute, the coordinator sends the power injection direction signals to the local controllers every 1 minute, and the local controllers communicate with their associated DER(s) every 10 seconds.

The scenario used to evaluate the functionality of the DERMS under different communication issues is chosen to be on April 3, 2019, from 10 a.m.--12 p.m. because it has a smooth PV and load profile that limits the exogenous effects compared to the implemented disturbance. 
At 11 a.m., we implement a step change disturbance in the VPP bound settings by changing the time-varying set point from 0.97 MW to 1.17 MW for Phase A, from 0.93 MW to 0.70 MW for Phase B, and from 0.95 MW to 0.61 MW for Phase C.
This is a plausible change that could be requested by the entity contracting the DERMS for some other higher level objective.
The evaluation metric described in Section \ref{sec:metrics} to determine the functionality of the DERMS is set with a time duration of 40 minutes, a voltage magnitude tolerance of 0.002 p.u., and a VPP tolerance of 0.02 MW.

\subsection{Upper Limit on Packet Losses}

We investigate the upper limit on the probability of a packet loss that keeps the DERMS functional for an individual communication category versus the percentage of communication channels with possible packet losses; all other channels have no possible packet loss.
Within HELICS, the random packet losses for the grid service measurements (voltage magnitudes and feeder head power measurements) are simulated by adding a native filter to the sending endpoint of the communication channels from the feeder to the DERMS coordinator; whereas the random packet losses for the power injection direction signals, sending PQ set points, and receiving PQ measurements are simulated within the local controller simulation modules.  

\begin{figure}
    \centering
    \includegraphics[width=0.95\columnwidth]{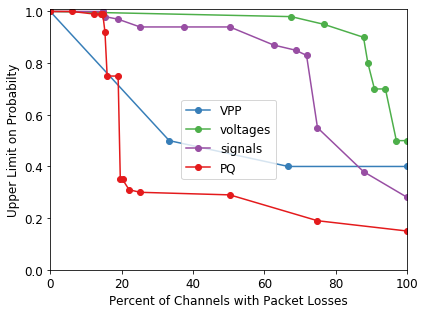}
    \caption{Upper limit on the probability of packet loss versus the percentage of channels with possible packet losses for each communication channel category.}
    \label{fig:Random_Drop_all}
\end{figure}

Figure \ref{fig:Random_Drop_all} gives the upper limit on the probability of a packet loss versus the percentage of channels with possible packet losses for each communication channel category.
The first observation is that when 100\% of the channels have a possible packet loss within their respective categories, the decreasing order with respect to the probability follows the order of the flow of information in the DERMS (see Figure \ref{fig:DERMS_diagram}).
The second observation is that as the percentage of channels with a possible packet loss increases, there is a threshold where the upper limit to the probability significantly drops, except for the VPP grid service, which has only three channels.
The decreasing order with respect to this threshold on the percentage of channels also follows the order of the flow of information in the DERMS.
This means that as the flow of the DERMS control information flows closer to the implementation of the power injections, it is less resistant to packet loss.

\subsection{Upper Limit on Link Failures}

We explore the upper limit on the probability of a link failure that keeps the DERMS functional versus the duration of the link failure.
Within HELICS, the random link failures are simulated within their specific simulation modules: The grid service measurement link failures are simulated within the feeder module; the power injection direction signal link failures are simulated within the DERMS coordinator module; and the local controller and DER link failures are simulated within the local controller modules.

\begin{figure}
    \centering
    \includegraphics[width=0.95\columnwidth]{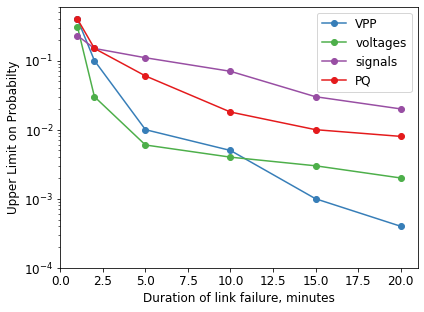}
    \caption{Upper limit on the probability of a link failure versus the duration of the failure for each communication channel category.}
    \label{fig:Random_Link_Failure_all}
\end{figure}

Figure \ref{fig:Random_Link_Failure_all} gives the upper limit on the probability of a link failure versus the duration of the failure for each communication channel category.
Both grid services (VPP, voltages) have a very fast drop from the interval from 1 to 5 minutes, and then they level off to a linear-logarithmic decrease, but the VPP has a steeper decline than the voltages.
This is likely because each phase of the feeder head powers has only a single point of failure, whereas each phase of the voltage magnitude measurements has many measurement nodes, and both are strongly correlated with each other in the same phase.
The reason that the decrease of the power injection signals is more gradual than all other categories is likely because each power injection direction signal carries a mix of information about the VPP and voltage magnitude bounds; thus, when some communication links are down, others still carry similar information.
The gradual leveling off of the PQ measurements and set point communication category is likely due to the DERs getting their set points stuck when a link fails, which are approximately where the set point would be if the link were operational.

\subsection{Upper Limit on Delays}

Finally, we investigate the upper bound on the delay added to each channel category.
Within HELICS, the communication delays for the grid service measurements are simulated by adding a native filter to the sending endpoint of the communication channels from the feeder to the DERMS coordinator, whereas the communication delays for the power injection direction signals, sending PQ set points, and receiving PQ measurements are simulated within the local controller modules.
The communication time step resolutions are changed for this communication issue compared to the others so that we can have more granular delay choices.
The grid service measurement communication and power injection direction signal communication are sent every 10 seconds, and communication between the local controllers and the DERs occurs every 2 seconds.

\begin{table}[]
    \centering
    \begin{tabular}{|l|c|}
        \hline
         \textbf{Communication Channel Category} & \textbf{Delay (seconds)} \\
         \hline\hline
         Feeder head power measurements & \\
         from the substation  & 80 \\
         to the DERMS coordinator & \\
         \hline
         Voltage magnitude measurements & \\
         from nodes on the feeder & 60 \\
         to the DERMS coordinator & \\
         \hline
         Power injection direction signals & \\
         from the DERMS coordinator & 50 \\
         to the DERMS local controllers &  \\
         \hline
         Power injection set point & \\
         from a DERMS local controller & 18 \\
         to its DERs &  \\
         \hline
         Power injection measurements &   \\
         from DERs &  $<0.2$ \\
         to their DERMS local controllers & \\
         \hline
    \end{tabular}
    \caption{Upper limit on the delay for each communication channel category.}
    \label{tab:delays}
\end{table}

Table \ref{tab:delays} gives the upper limit on the delay for each communication channel category.
Just like in the upper limit on the probability of a packet loss, the decreasing order with respect to the delay size follows the order of the flow of information in the DERMS, which means that as the flow of the DERMS control information gets closer to the implementation of the power injections, it is less resistant to delays.
In fact, we could not find a size of the delay at or exceeding 200 ms for the communication of the power injection measurements from the DER to the local controller that keeps the DERMS functional due to the 200 ms of simulation time that the simulator uses to pass messages from the OpenDSS module to the local controllers (see Figure \ref{fig:HELICS_diagram}).
With this observation in mind, we suggest that the DERs participating in the DERMS should be continuously sending their power injection measurements to their DERMS local controllers at a faster frequency than the local controller makes set point decisions so that the local controller can confirm the true power injection before deciding a new set point.

\section{Conclusion}
This paper studies the upper limits on the severity of issues that can occur in the communication systems used by a primal-dual-based DERMS, which include packet losses, link failures, and delays.
This information can be used by a DERMS operator to decide on the necessary specifications for the communication infrastructure it invests to have a functional DERMS operation.
We found the upper limits through numerical simulations of a primal-dual-based DERMS providing grid services to a distribution feeder with approximately 2,000 nodes by controlling 163 PV generators and 140 energy storage batteries.
In general, we found that the DERMS becomes less resistant to communication issues as the information flows closer to the DERs implementing the controllable power injections.
Fortunately, this means that specifications are more restrictive for the HANs communicating with the DERs than the more complex NAN needed to communicate with the distribution feeder and the DERMS operator.

\section*{Acknowledgments}
This work was authored by the National Renewable Energy Laboratory, operated by Alliance for Sustainable Energy, LLC, for the U.S. Department of Energy (DOE) under Contract No. DE-AC36-08GO28308. Funding provided by U.S. Department of Energy Office of Energy Efficiency and Renewable Energy Solar Energy Technologies Office Award Number TCF-21-25008. The views expressed in the article do not necessarily represent the views of the DOE or the U.S. Government. The U.S. Government retains and the publisher, by accepting the article for publication, acknowledges that the U.S. Government retains a nonexclusive, paid-up, irrevocable, worldwide license to publish or reproduce the published form of this work, or allow others to do so, for U.S. Government purposes.


\end{document}